\documentclass[conference]{IEEEtran}
\IEEEoverridecommandlockouts
\usepackage{cite}
\usepackage{amsmath,amssymb,amsfonts}
\usepackage{algorithm,algorithmic}
\usepackage[dvipdfm]{graphicx}
\usepackage{textcomp}
\usepackage{xcolor}
\def\BibTeX{{\rm B\kern-.05em{\sc i\kern-.025em b}\kern-.08em
    T\kern-.1667em\lower.7ex\hbox{E}\kern-.125emX}}
\begin{document}
%-----------------------------------%
% Title
%-----------------------------------%
\title{Acceleration of complex matrix multiplication using arbitrary precision floating-point arithmetic
%\thanks{This study was supported by JSPS KAKENHI Grant Number 23K11127.}
}

\author{\IEEEauthorblockN{Tomonori Kouya}
%\author{Anonymous}
\IEEEauthorblockA{\textit{Shizuoka Institute of Science and Technology} \\
Fukuroi, Japan \\
ORCID: 0000-0003-0178-5519}
}

\maketitle

%-----------------------------------%
% Abstract
%-----------------------------------%
\begin{abstract}
Efficient multiple precision linear numerical computation libraries such as MPLAPACK are critical in dealing with ill-conditioned problems. Specifically, there are optimization methods for matrix multiplication, such as the Strassen algorithm and the Ozaki scheme, which can be used to speed up computation. For complex matrix multiplication, the 3M method can also be used, which requires only three multiplications of real matrices, instead of the 4M method, which requires four multiplications of real matrices. In this study, we extend these optimization methods to arbitrary precision complex matrix multiplication and verify the possible increase in computation speed through benchmark tests. The optimization methods are also applied to complex LU decomposition using matrix multiplication to demonstrate that the Ozaki scheme can be used to achieve higher computation speeds.
\end{abstract}

%-----------------------------------%
% Keywords
%-----------------------------------%
\begin{IEEEkeywords}
complex matrix multiplication, arbitrary precision arithmetic, Strassen algorithm, Ozaki scheme
\end{IEEEkeywords}

%-----------------------------------%
% Introduction
%-----------------------------------%
\section{Introduction}
Efficient multiple precision (MP) floating-point arithmetic is critical for a linear computation library, such as MPLAPACK~\cite{mplapack}, to obtain accurate numerical solutions of ill-conditioned problems. Currently, MPLAPACK is the de-facto standard MP linear library owing to embedded trustworthy MP floating-point arithmetic libraries such as QD~\cite{qd}, GNU MP~\cite{gmp}, MPFR~\cite{mpfr}, and MPC~\cite{mpc}.  Nevertheless, MPBLAS of MPLAPACK, which is based on Reference BLAS (Basic Linear Algebra Subprograms), has plenty of room for improvement. The current LAPACK using IEEE binary32 and binary64 is utilized in conjunction with highly optimized BLAS such as ATLAS~\cite{atlas}, OpenBLAS~\cite{openblas}, and Intel Math Kernel~\cite{imkl}. The optimized MPBLAS library can accelerate the functionalities of MPLAPACK even more. Therefore we can expect that this type of optimization when applied to complex basic linear computations is more effective than the original MPBLAS.

In speeding up multiple precision linear computations using the MPFR~\cite{mpfr} arbitrary precision floating-point library, which is based on the multiple precision natural number (MPN) kernel provided by GNU MP, the algorithm that reduces the amount of computation is the most effective in speeding up the process. For matrix multiplication, there are optimization methods such as the Strassen algorithm \cite{strassen_original} and the Ozaki scheme~\cite{ozaki_scheme}, which can be used to speed up calculations. In particular, the recently developed Ozaki scheme is an algorithm that pursues both accuracy and speed by dividing the original matrix into matrices of short mantissa parts, thereby performing fast low-precision matrix multiplication without errors. The Ozaki scheme has already been shown to be effective in float128 calculations~\cite{mukunoki_binary128}; however, since its effectiveness depends on the nature of the matrices used, we must verify the effectiveness of the matrices for concrete problems through benchmark tests. In addition, although there is ongoing research on optimization of real matrix multiplication, there are no comparative studies on the effectiveness of these optimization techniques for multiple precision complex matrix multiplication.

The 3M method is used in MPC~\cite{mpc}, an arbitrary precision complex arithmetic library based on MPFR, which can be employed in complex matrix multiplication. It is expected that further speed-up can be achieved by using the 3M method. There are already comparisons of BLIS implementations for binary32 and binary64 complex matrix multiplications\cite{cmatmul_3m4m}, and it has been shown that the 3M method is effective in reducing computational time despite its numerical instability. Another example of implementing binary64 complex multiplication based on the Ozaki scheme for real matrices is that of Kazal.et.al~\cite{cmatmul_ozaki}, which uses the 4M method and, hence, cannot be compared with the 3M method.

In this study, we implement complex matrix multiplication using the Strassen algorithm and the Ozaki scheme based on the 3M and 4M methods and report the results of benchmark tests on well-conditioned matrices generated using random numbers. We have already verified the speedup of arbitrary precision real matrix multiplication based on the Strassen algorithm and Ozaki scheme using MPFR\cite{kouya_utsugiri_ozaki}, and we expect to be able to speed up the 3M and 4M methods based on these results so far. As an application, we implemented a complex LU decomposition incorporating these matrix multiplications and compared the computational time and accuracy of the numerical solutions through benchmark tests on well-conditioned linear system of equations. As a result, it is found that the Ozaki scheme is effective for relatively small precision; however, the computational time increases as precision increases and is slower than the Strassen matrix multiplication for precisions greater than 768 bits.

The following two computation environments, EPYC and Xeon, are used in this study. MPLAPACK and our library, including QD and GNU MP, MPFR, and MPC, are natively compiled with Intel Compiler and DGEMM (cblas\_dgemm function) in Intel Math Kernel. 
\begin{description}
	\item[Xeon] Intel Xeon W-2295 3.0GHz 18 cores, Ubuntu 20.04.3 LTS, Intel Compiler version 2021.5.0, MPLAPACK 2.0.1, GNU MP 6.2.1, MPFR 4.1.0, MPC 1.2.1
	\item[EPYC] AMD EPYC 7402P 24 cores, Ubuntu 18.04.6 LTS, Intel Compiler version 2021.4.0, MPLAPACK 2.0.1, GNU MP 6.2.1, MPFR 4.1.0, MPC 1.2.1
\end{description}

% ---------------------------------------
\section{Mathematical notation}
% ---------------------------------------
Here we define $\mathbb{F}_{bS}$ and $\mathbb{F}_{bL}$ as sets of the $S$- and $L$-bit mantissas of floating-point numbers, respectively. For instance, $\mathbb{F}_{b24}$ and $\mathbb{F}_{b53}$ refer to sets of IEEE754-1985 binary32 and binary64 floating-point numbers, respectively. Although any mantissa length can be selected in MPFR arithmetic, the set of MPFR numbers is expressed as $\mathbb{F}_{bM}$, which is primarily defined as $M$-bit using the mpfr\_set\_default\_prec function.

We used $(\mathbf{x})_i$ $(= x_i)$ as the $i$-th element of the $n$-dimensional real or complex vector $\mathbf{x} = [x_i]_{i=1, 2, .., n}\in\mathbb{R}^n$ or $\mathbb{C}^n$, and real or complex matrix $(A)_{ij} (= a_{ij})$ as the $(i, j)$-th element of $A = [a_{ij}]_{i=1, 2, ..., m, j = 1, 2, ..., n}\in\mathbb{R}^{m\times n}$ or $\mathbb{C}^{m\times n}$.

Real and imaginary parts of complex vector $\mathbf{x}\in\mathbb{C}^n$ and complex matrix $A\in\mathbb{C}^{m\times n}$ are expressed as $\mathrm{Re}(\mathbf{x}), \mathrm{Im}(\mathbf{x})\in\mathbb{R}^n$, and $\mathrm{Re}(A), \mathrm{Im}(A)\in\mathbb{R}^{m\times n}$, respectively.

%-----------------------------------%
% 
%-----------------------------------%
\section{Algorithms of complex matrix multiplication}

In this section, the algorithms used in this paper are introduced. First, the 4M and 3M methods of complex multiplication are explained and subsequently, their application to matrix multiplication is described. Next, the Strassen matrix multiplication and the Ozaki scheme are explained.

\subsection{3M and 4M methods}

Given $A\in\mathbb{C}^{m\times l}$ and $B\in\mathbb{C}^{l\times n}$, we calculate the complex product 
\begin{equation}
	C := AB\in\mathbb{C}^{m\times n}. \label{eqn:cmatmul}
\end{equation}
From the mathematical definition of the above matrix multiplication, we can obtain the element of $C$ as follows:
\begin{equation}
	(C)_{ij} := \sum_{k=1}^l (A)_{ik} (B)_{kj}, \label{eqn:simple_triple_loop}
\end{equation}
where the element $(C)_{ij}$ can be calculated with a simple triple-loop in the programs. The formula (\ref{eqn:simple_triple_loop}) includes complex multiplication and addition. The complex multiplication is calculated as follows:
\begin{equation}
	\begin{split}
		\mathrm{Re}((A)_{ik} (B)_{kj}) &:= \mathrm{Re}((A)_{ik})\mathrm{Re}((B)_{kj}) \\
			& - \mathrm{Im}((A)_{ik})\mathrm{Im}((B)_{kj}) \\
		\mathrm{Im}((A)_{ik} (B)_{kj}) &:= \mathrm{Re}((A)_{ik})\mathrm{Im}((B)_{kj}) \\
			& + \mathrm{Im}((A)_{ik})\mathrm{Re}((B)_{kj})
	\end{split} \label{eqn:4m_way_base}
\end{equation}
The formula for standard complex multiplication is called the ``4M" method.

In contrast, the method shown below, which reduces the number of multiplications by one, as in the Karatsuba method, is called the ``3M" method.
\begin{equation}
	\begin{split}
		t_1 &:= \mathrm{Re}((A)_{ik})\mathrm{Re}((B)_{kj}) \\
		t_2 &:= \mathrm{Im}((A)_{ik})\mathrm{Im}((B)_{kj}) \\
		\mathrm{Re}((A)_{ik} (B)_{kj}) &:= t_1 - t_2 \\
		\mathrm{Im}((A)_{ik} (B)_{kj}) &:= (\mathrm{Re}((A)_{ik}) + \mathrm{Im}((A)_{ik})) \\
		& \cdot (\mathrm{Re}((B)_{kj}) + \mathrm{Im}((B)_{kj})) \\
		& - t_1 - t_2
	\end{split} \label{eqn:3m_way_base}
\end{equation}
In the MPC library based on MPFR arithmetic, the fast complex multiplication (mpc\_mul function) is implemented using the 3M method. All multiple precision complex linear computations defined in MPLAPACK/MPBLAS are constructed on ``mpcomplex" C++ class with MPC arithmetic functions.

The 3M and 4M complex general matrix multiplication (CGEMM) methods are implemented using real matrix multiplication. The 4M CGEMM method is expressed in a straightforward manner, as shown in (\ref{eqn:4m_way_base}), as follows:
\begin{equation}
	\begin{split}
		\mathrm{Re}(AB) &:= \mathrm{Re}(A)\mathrm{Re}(B) - \mathrm{Im}(A)\mathrm{Im}(B) \\
		\mathrm{Im}(AB) &:= \mathrm{Re}(A)\mathrm{Im}(B) + \mathrm{Im}(A)\mathrm{Re}(B)
	\end{split}\ \label{eqn:4m_way}
\end{equation}
The 3M CGEMM method is also simple to express, as shown in (\ref{eqn:3m_way_base}), as follows:
\begin{equation}
	\begin{split}
		T_1 &:= \mathrm{Re}(A)\mathrm{Re}(B) \\
		T_2 &:= \mathrm{Im}(A)\mathrm{Im}(B) \\
		\mathrm{Re}(AB) &:= T_1 - T_2 \\
		\mathrm{Im}(AB) &:= (\mathrm{Re}(A) + \mathrm{Im}(A))(\mathrm{Re}(B) + \mathrm{Im}(B)) \\
		&- T_1 - T_2
	\end{split} \label{eqn:3m_way}
\end{equation}
\tablename\ \ref{table:num_calc_3m4m} describes the number of real arithmetic operations to obtain complex number or matrix products using the 3M and 4M methods. For complex multiple precision arithmetic and matrix multiplication, reducing one multiplication is more efficient than increasing three additions/subtractions. Although it is possible to obtain errors in imaginary parts~\cite{higham_accuracy}, we could not confirm this phenomenon in our benchmark tests, which are described later.

\begin{table}
	\begin{center}
		\caption{The number of real calculations in 4M and 3M methods}\label{table:num_calc_3m4m}
		\begin{tabular}{|c|c|c|}\hline
			& ADD, SUB & MUL \\ \hline
		4M  & 2 & 4 \\ \hline
		3M  & 5 & 3 \\ \hline
		\end{tabular}
	\end{center}
\end{table}

%\subsection{4M and 3M ways}

\subsection{Strassen algorithm}

The Strassen matrix multiplication algorithm\cite{strassen_original} is categorized as a divide-and-conquer method and is well-known to drastically reduce the number of arithmetic operations when the size of matrices is increased. Multiple precision floating-point arithmetic incurs high costs, so the Strassen algorithm shown in Algorithm \ref{algo:strassen} is useful when multiple precision matrix multiplication is needed. We have implemented the Strassen algorithm using the MPC library. The thresholds of matrix size are defined as $m_0 = n_0 = 32$.

%\begin{figure}[htb]
	\begin{algorithm}[htb]
		\algsetup{linenosize=\small}
		\small
		\caption{Strassen algorithm for complex matrix multiplication}\label{algo:strassen}
		\textbf{Input:} $A = [A_{ij}]_{i,j=1,2}\in\mathbb{C}^{m\times l}$, $A_{ij}\in\mathbb{C}^{m/2\times l/2}$, $B = [B_{ij}]_{i,j=1,2}\in\mathbb{C}^{l\times n}$, $B_{ij}\in\mathbb{C}^{l/2\times n/2}$ \\
		\textbf{Output:} $C:=$\mbox{Strassen}$(A, B)=AB\in\mathbb{C}^{m\times n}$
		\begin{algorithmic}
			\IF{$m < m_0$ \&\& $n < n_0$}
				\STATE $C := AB$
			\ENDIF
			\STATE $P_1 := \mbox{Strassen}(A_{11} + A_{22}, B_{11} + B_{22})$
			\STATE $P_2 := \mbox{Strassen}(A_{21} + A_{22}, B_{11})$
			\STATE $P_3 := \mbox{Strassen}(A_{11}, B_{12} - B_{22})$
			\STATE $P_4 := \mbox{Strassen}(A_{22}, B_{21} - B_{11})$
			\STATE $P_5 := \mbox{Strassen}(A_{11} + A_{12}, B_{22})$
			\STATE $P_6 := \mbox{Strassen}(A_{21} - A_{11}, B_{11} + B_{12})$
			\STATE $P_7 := \mbox{Strassen}(A_{12} - A_{22}, B_{21} + B_{22})$
			\STATE $C_{11} := P_1 + P_4 - P_5 + P_7$; $C_{12} := P_3 + P_5$
			\STATE $C_{21} := P_2 + P_4$; $C_{22} := P_1 + P_3 - P_2 + P_6$
			\STATE $C := [C_{ij}]_{i,j=1,2}$
		\end{algorithmic}
	\end{algorithm}
	%\end{figure}

\subsection{Ozaki scheme}

The usefulness of the Ozaki scheme\cite{ozaki_scheme} is also becoming clear at multiple precision levels owing to the success of Mukunoki et al. in accelerating the float128 precision matrix multiplication\cite{mukunoki_binary128}. %The float128 arithmetic, which is supported by GCC, exhibits TD to QD precision performance for addition and multiplication, and is expected to be appropriate for precision levels within this neighborhood. 
The float128 arithmetic supported by GCC, features triple-double to quadruple-double precision performance for addition and multiplication, which is expected to be sufficient for this precision range.

The Ozaki scheme is an algorithm that aims to simultaneously accelerate performance and improve accuracy by dividing original matrices into another matrices with elements represented by shorter digits. Similar to the ``Split'' method in error-free transformation technique, this approach takes advantage of the speed of optimized short-precision matrix multiplication (xGEMM) functions without round-off errors. For a given matrix $A \in \mathbb{R}^{m\times l}$ and $B \in \mathbb{R}^{l\times n}$, to obtain a matrix product $C := AB \in \mathbb{R}^{m\times n}$ of long $L$-bit precision, $A$ and $B$ are divided using the Ozaki scheme, where $d \in \mathbb{N}$ is the maximum number of divisions of short $S$-bit precision matrices ($S << L$), as shown in Algorithm \ref{algo:ozaki_scheme}. The $S$-bit arithmetic is used for calculations when no particular description is provided, and the $L$-bit arithmetic is used only when high-precision operations are required.

%\begin{figure}[htb]
\begin{algorithm}[htb]
    \algsetup{linenosize=\small}
    \small
    \caption{Ozaki scheme for multiple precision matrix multiplication}\label{algo:ozaki_scheme}
	\textbf{Input:} $A \in \mathbb{F}_{bL}^{m\times l}, B \in\mathbb{F}_{bL}^{l\times n}$ \\
	\textbf{Output:} $C \in \mathbb{F}_{bL}^{m\times n}$
    \begin{algorithmic}
        \STATE $A^{(S)} := A$, $B^{(S)} := B$ : $A^{(S)} \in \mathbb{F}_{bS}^{m\times l}$, $B^{(S)} \in \mathbb{F}_{bS}^{l\times n}$
		% divide A and B
		\STATE $\mathbf{e} := [1\ 1\ ...\ 1]^T\in \mathbb{F}_{bS}^l$
		\STATE $\alpha := 1$
		\WHILE{$\alpha < d$}
			\STATE ${\boldsymbol\mu}_A := [\max_{1 \leq p \leq l} |(A^{(S)})_{ip}|]_{i=1, 2, ..., m} \in \mathbb{F}_{bS}^m$
			\STATE ${\boldsymbol\mu}_B := [\max_{1 \leq q \leq l} |(B^{(S)})_{qj}|]_{j=1, 2, ..., n} \in \mathbb{F}_{bS}^n$
			\STATE ${\boldsymbol\tau}_A := [2^{\lceil \log_2(({\boldsymbol\mu}_A)_i)\rceil + \lceil (S + \log_2(l)) / 2 \rceil} ]_{i = 1, 2, ..., m} \in \mathbb{F}_{bS}^m$
			\STATE ${\boldsymbol\tau}_B := [2^{\lceil \log_2(({\boldsymbol\mu}_B)_j)\rceil + \lceil (S + \log_2(l)) / 2 \rceil} ]_{j = 1, 2, ..., n} \in \mathbb{F}_{bS}^n$
			\STATE $S_A := \boldsymbol\tau_A \mathbf{e}^T$
			\STATE $S_B := \mathbf{e} \boldsymbol\tau_B^T$
			\STATE $A_\alpha := (A^{(S)} + S_A) - S_A$: $A_\alpha \in \mathbb{F}_{bS}^{m\times l}$
			\STATE $B_\alpha := (B^{(S)} + S_B) - S_B$: $B_\alpha \in \mathbb{F}_{bS}^{l\times n}$
			\STATE $A := A - A_\alpha$, $B := B - B_\alpha$ : $L$-bit FP arithmetic
			\STATE $A^{(S)} := A$, $B^{(S)} := B$
			\STATE $\alpha := \alpha + 1$  
		\ENDWHILE
		% multiplication
		\STATE $A_d := A^{(S)}$, $B_d := B^{(S)}$
		\STATE $C := O$
		\FOR{$\alpha = 1, 2, ..., d$}
			\FOR{$\beta = 1, 2, ..., d - \alpha + 1$}
				\STATE $C_{\alpha\beta} := A_\alpha B_\beta$
			\ENDFOR
			\STATE $C := C + \sum^{d - \alpha + 1}_{\beta = 1} C_{\alpha\beta}$ : $L$-bit FP arithmetic
		\ENDFOR
    \end{algorithmic}
\end{algorithm}
%\end{figure}

We implemented complex matrix multiplication based on 4M CGEMM (\ref{eqn:4m_way}) and 3M CGEMM (\ref{eqn:3m_way}) methods including the Ozaki scheme to obtain real matrix products. For practical purposes, it is desirable to set the maximum number of divisions $d$ for both the real and imaginary parts of the calculation; however, in this case, the calculations were performed on a real example using values of the same order in both cases, and a common $d$ is used.

%-----------------------------------%
% 
%-----------------------------------%
\section{Benchmark test}

In this section, we present the results of complex matrix multiplication in solving complex linear system of equations including a complex coefficient matrix to be well-conditioned generated using random numbers. In particular, since the accuracy of the Ozaki scheme depends on the number of partitions, it is necessary to compare the computational time for the smallest number of partitions $d$ to obtain the best accuracy. For this reason, graphs and tables presenting both computational time and computation accuracy are included.

%-----------------------------------%
% 
%-----------------------------------%
\subsection{Complex matrix multiplication}

Here we describe a benchmark test using complex square matrix multiplication. The elements of $A$ and $B$ used in the matrix multiplication (\ref{eqn:cmatmul}) are random numbers (using the mpfr\_nrandom function) with the real and imaginary parts following a normal distribution of $[-1, 1]$. The matrix size is increased and the computational time and maximum relative error per element of $C$ are compared. The sorts of precision used are 256, 512, and 768 bits. For comparison, the results of the Cgemm function of MPBLAS are also included, omitting cases with excessive computational time.

\begin{figure}%[htb]
    \begin{center}
        \includegraphics[width=.225\textwidth]{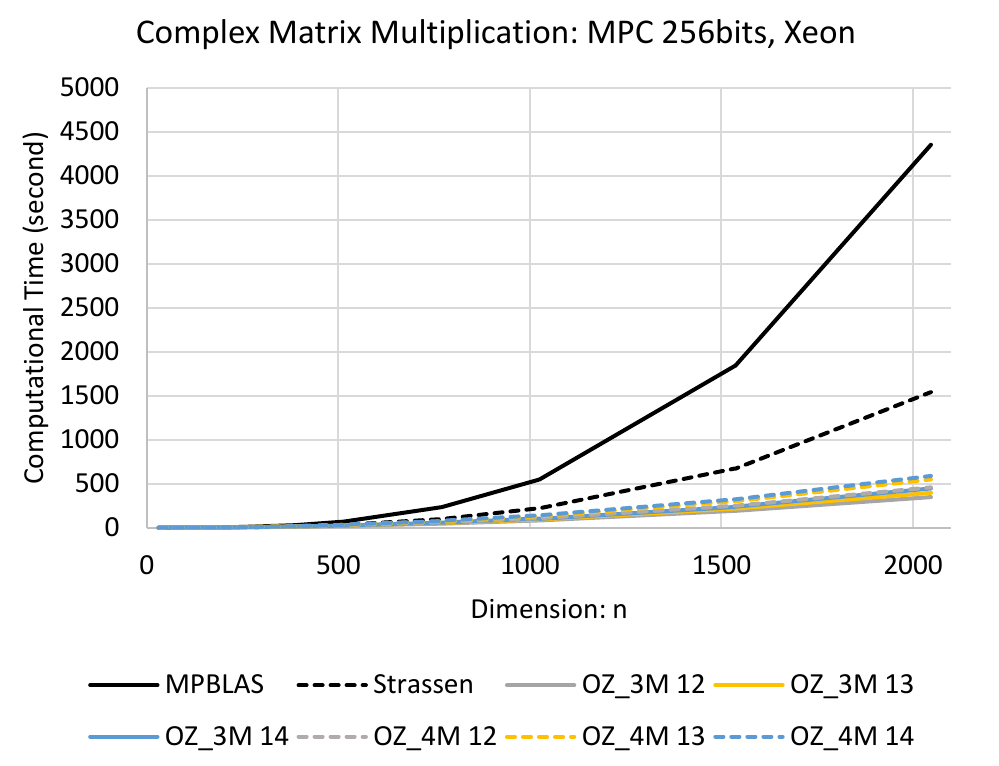}
        \includegraphics[width=.225\textwidth]{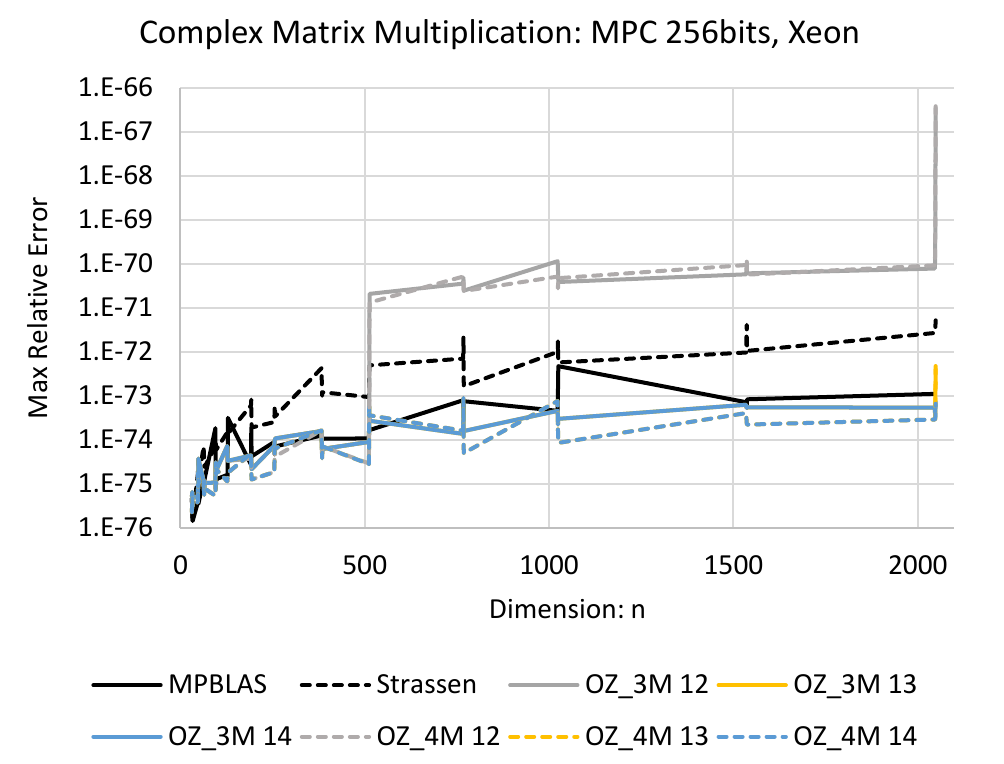}

		\includegraphics[width=.225\textwidth]{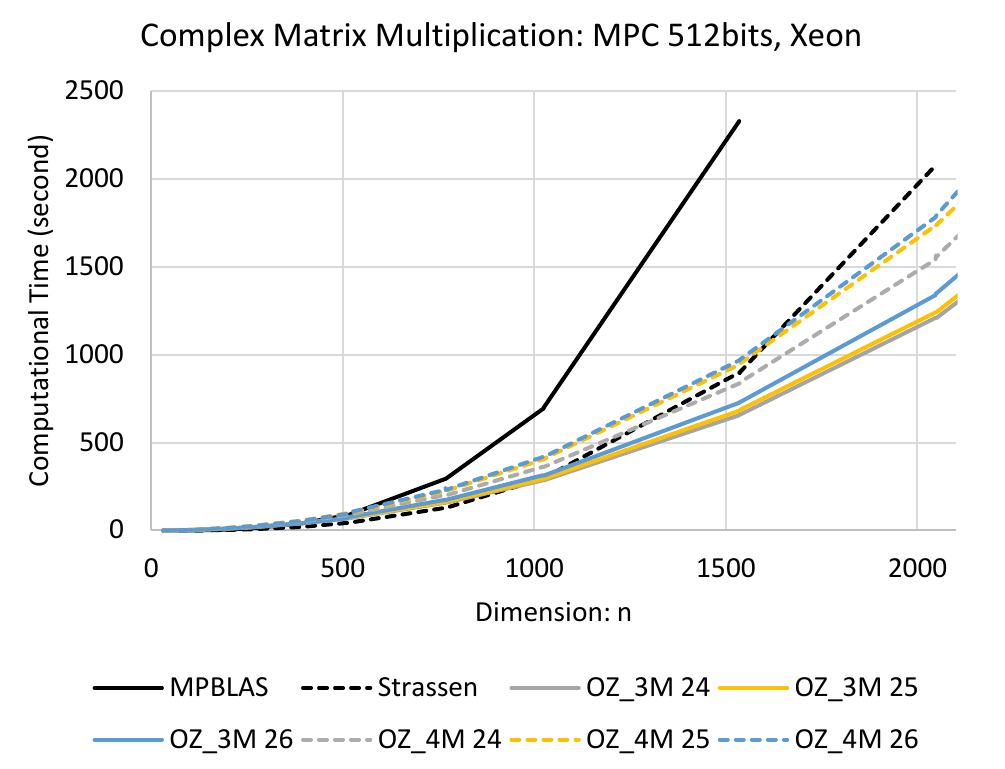}
        \includegraphics[width=.225\textwidth]{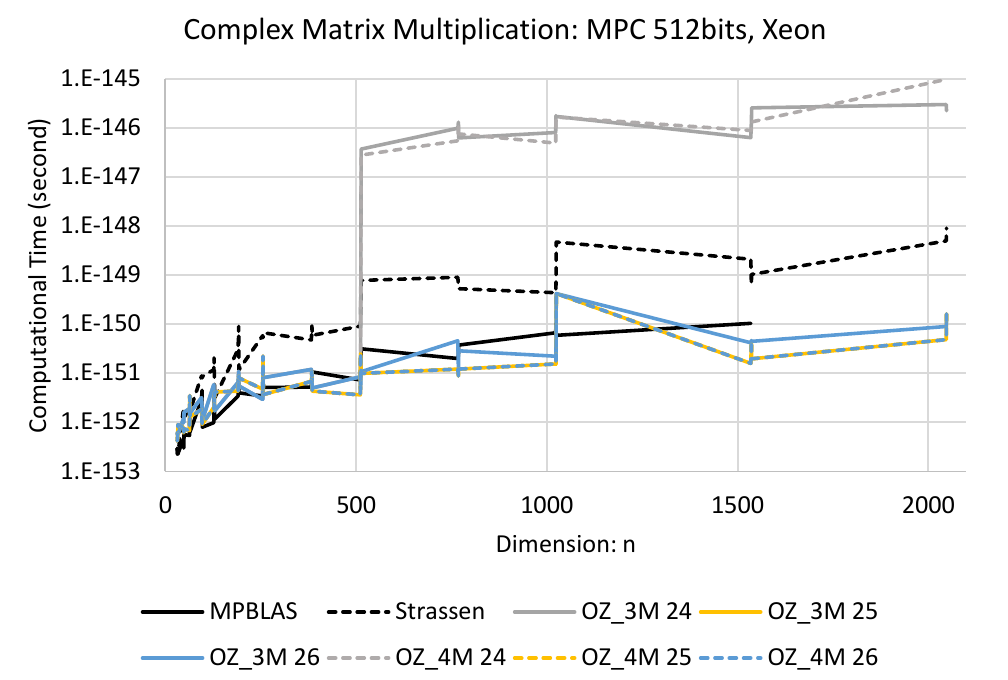}

		\includegraphics[width=.225\textwidth]{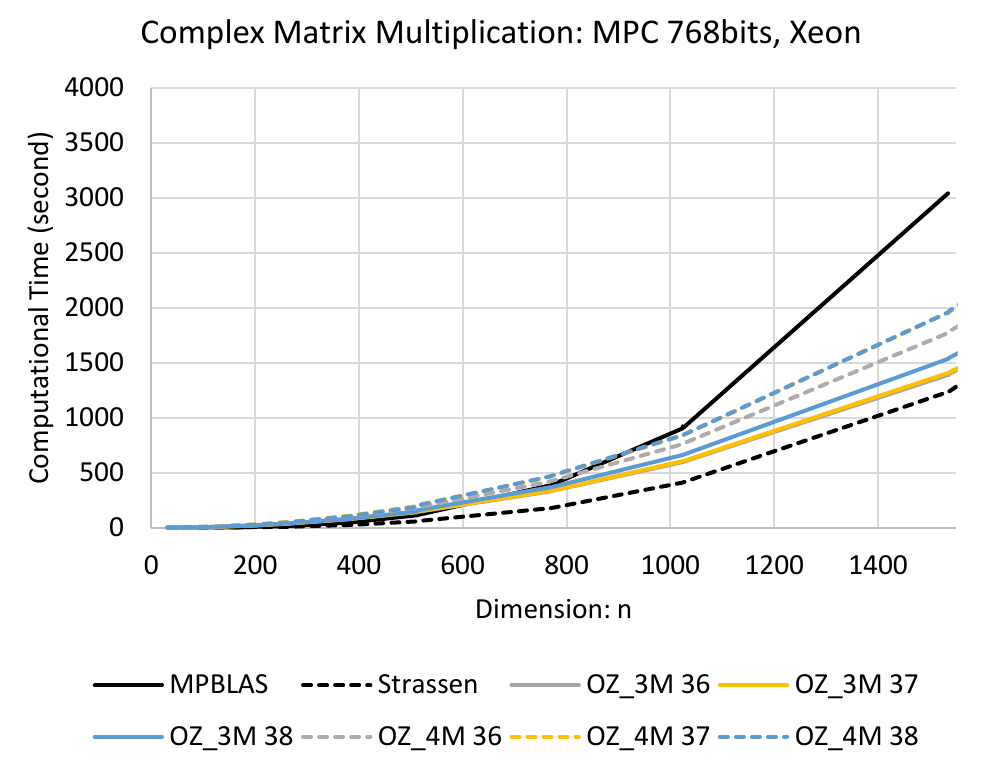}
        \includegraphics[width=.225\textwidth]{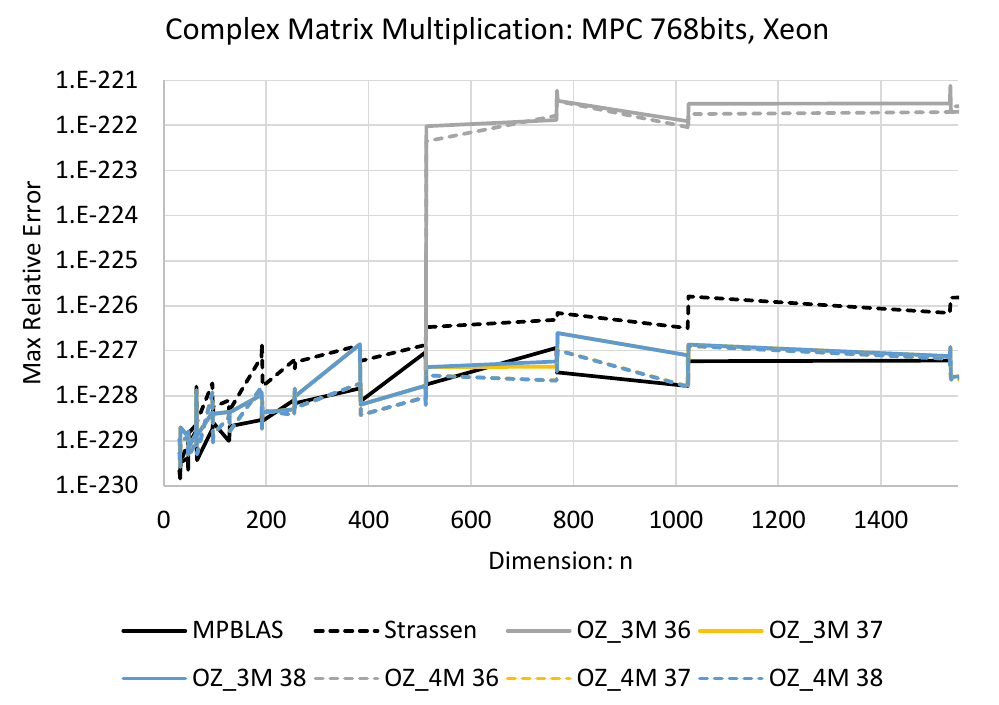}

		\caption{Computational time (left) and maximum relative error (right) of complex matrix multiplication on Xeon}
        \label{fig:matmulbench_mpc_xeon}
    \end{center}
\end{figure}

The graphs of computational time (left) and maximum relative error (right) on the Xeon environment are shown in \figurename\ref{fig:matmulbench_mpc_xeon}. In these figures, ``OZ\_4M" and ``OZ\_3M" mean results obtained by the Ozaki scheme based on 4M and 3M methods, respectively, and the numbers that follow mean the maximum numbers of divisions, $d$.

In the graphs of maximum relative error, the errors increase as matrix size increases. Overall, the Strassen matrix multiplication has the largest error, and the accuracy deteriorates by up to two decimal places compared to that of other CGEMM. 

In the Ozaki scheme, we also see that $d \geq 13$ for 256-bit, $d \geq 25$ for 512-bit, and $d \geq 37$ for 768-bit calculations minimize the relative error for any matrix size. Although this example did not require a large number of divisions for the Ozaki scheme because of the small order difference in the absolute values of the matrix elements, for DGEMM, a larger number of divisions is still required with increasing precision. 

We cannot confirm that a relative error of more than one decimal place occurred between the 3M and 4M CGEMM methods. The 3M CGEMM method is faster than the 4M CGEMM method by a 24 to 26\% decrease in computational time.

\begin{figure}%[htb]
    \begin{center}
        \includegraphics[width=.225\textwidth]{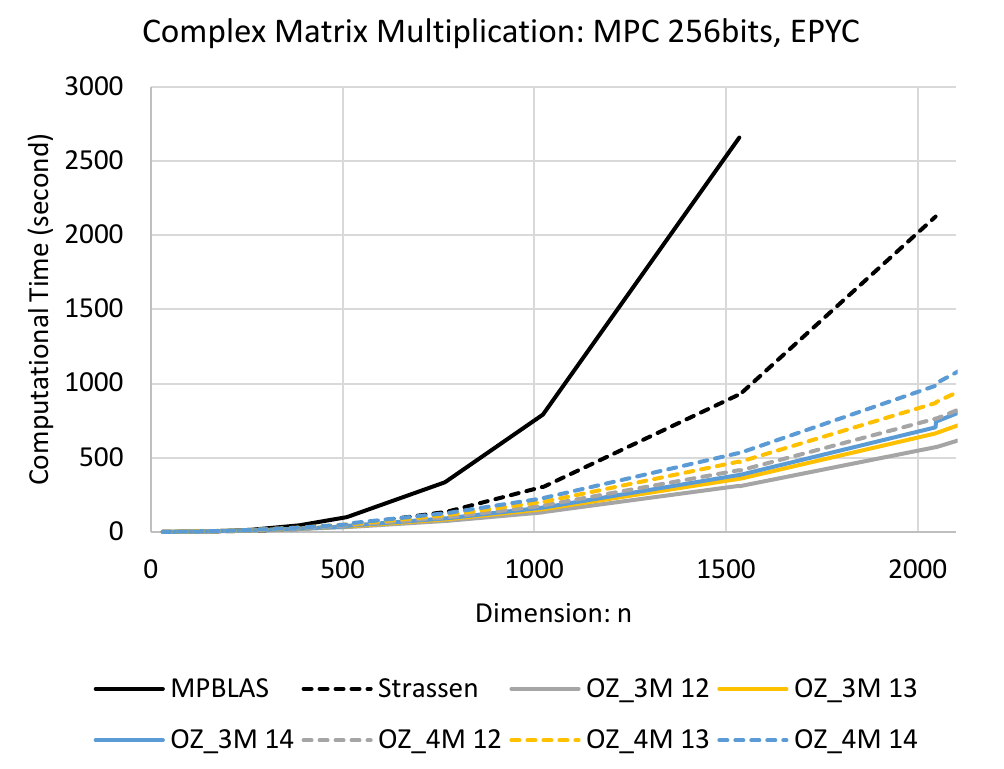}
		\includegraphics[width=.225\textwidth]{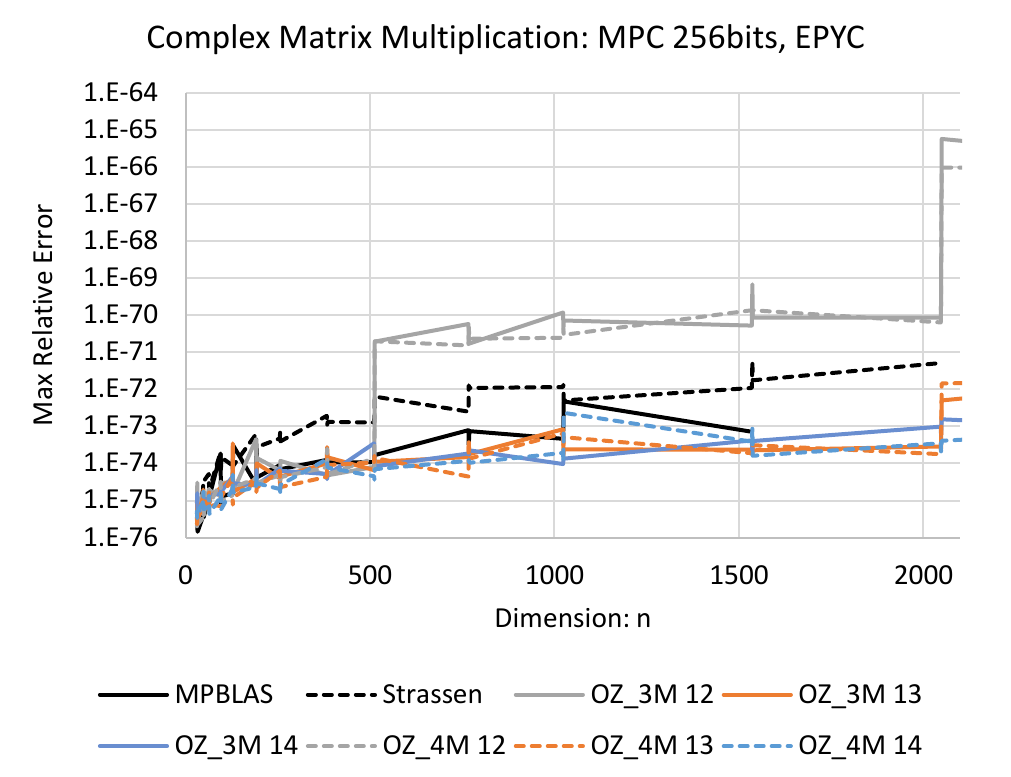}

		\includegraphics[width=.225\textwidth]{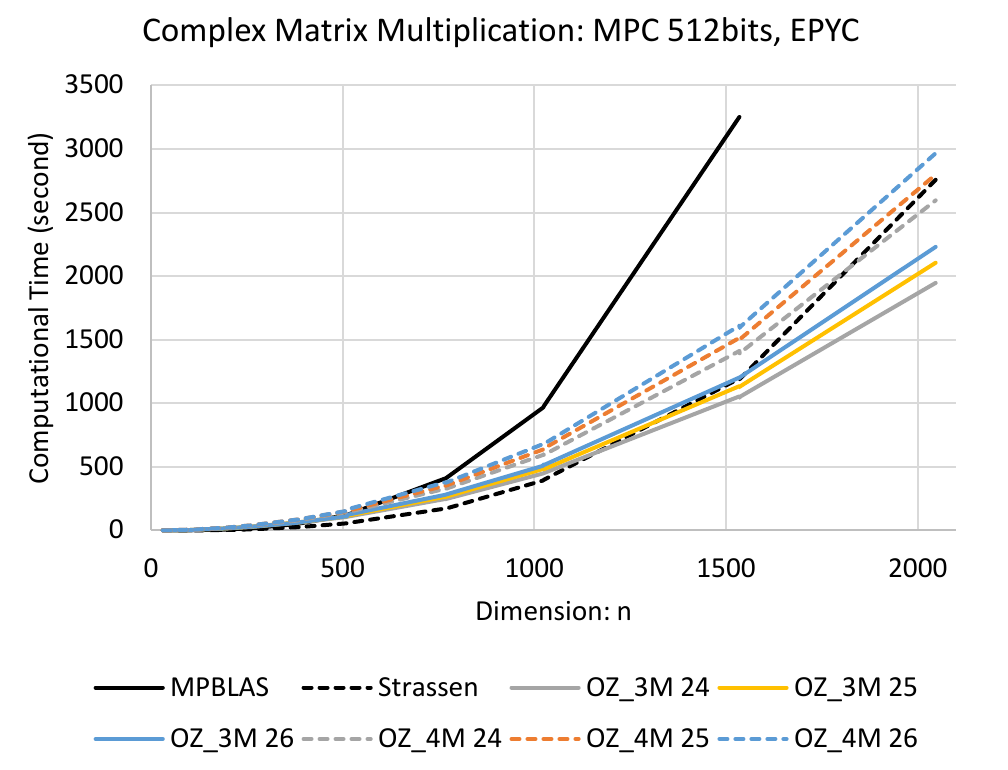}
		\includegraphics[width=.225\textwidth]{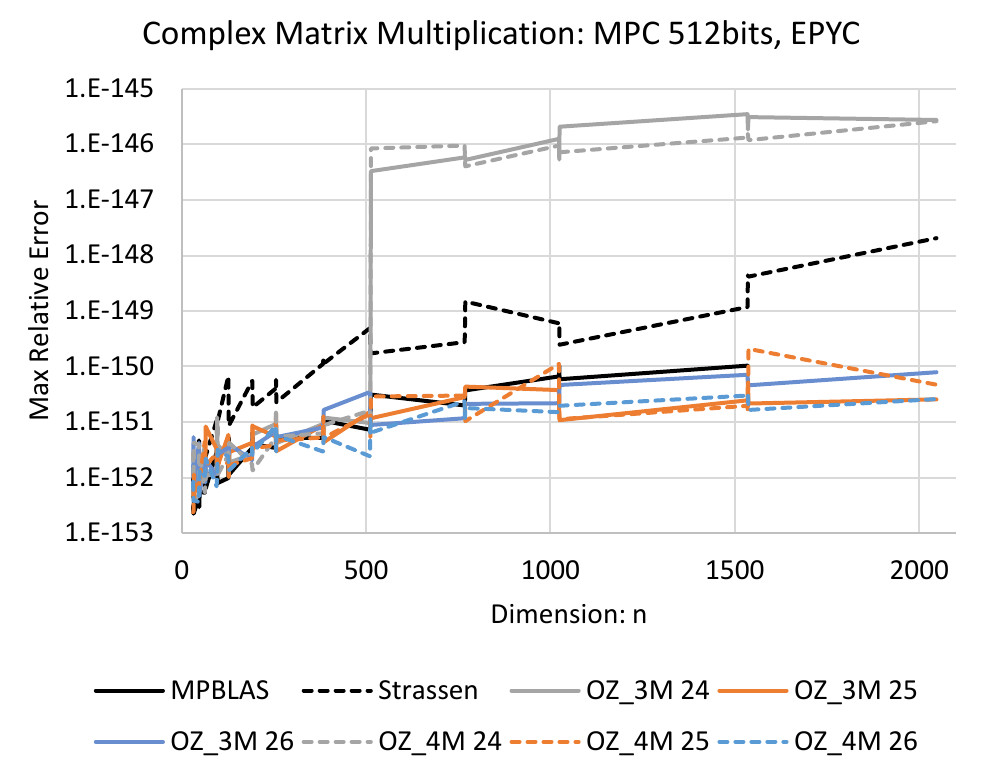}

		\includegraphics[width=.225\textwidth]{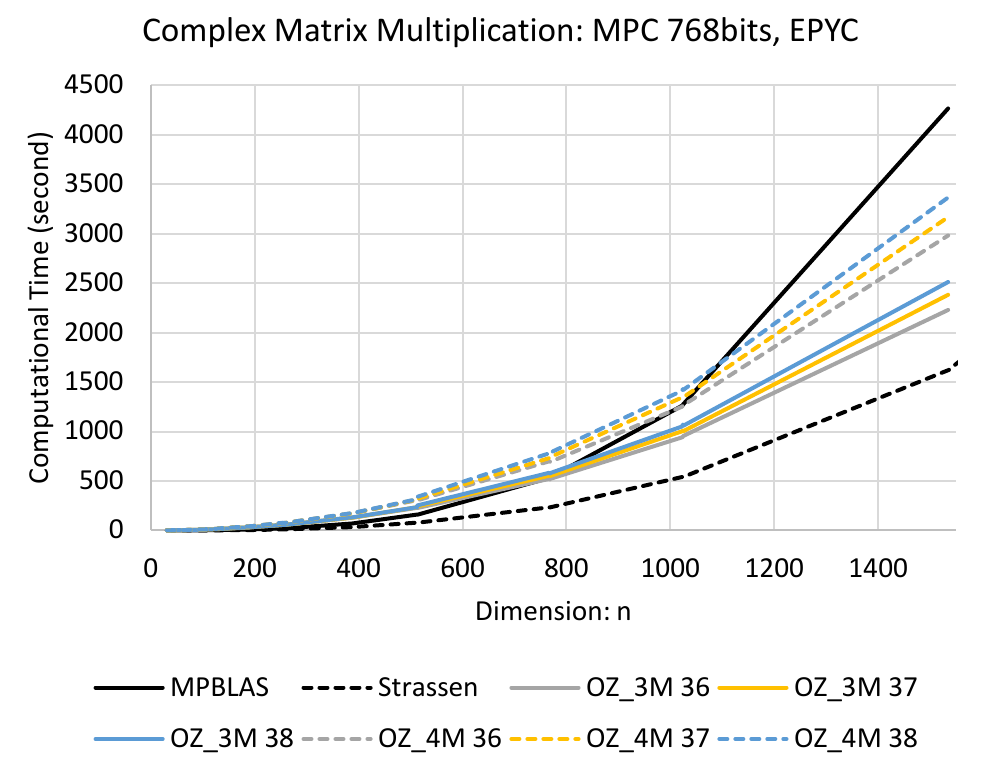}
		\includegraphics[width=.225\textwidth]{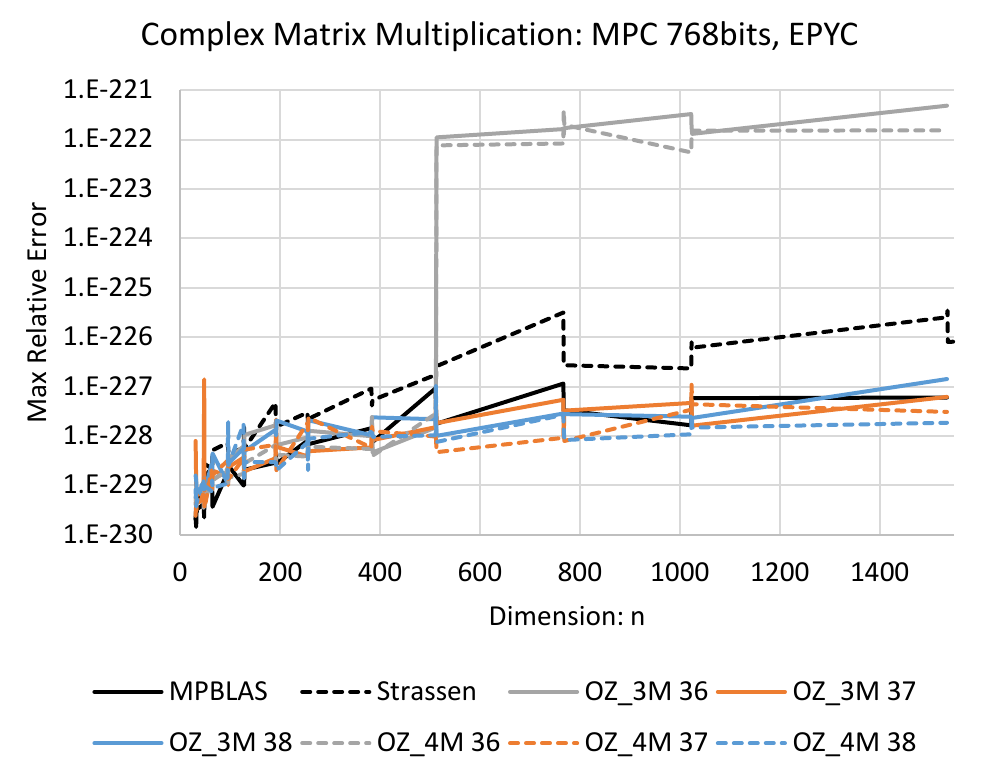}

		\caption{Computational time (left) and relative error (right) of complex matrix multiplication on EPYC}
        \label{fig:matmulbench_mpc_epyc}
    \end{center}
\end{figure}

Next, the graphs of computational time (left) and maximum relative error (right) are shown in \figurename\ref{fig:matmulbench_mpc_epyc} for the EPYC environment.

Overall, the computational time is slower than on the Xeon: about 20\% slower for MPC alone, and up to 40\% slower for the Ozaki scheme with DGEMM. There is little difference in the maximum relative error. In terms of computational time, the Ozaki scheme is the fastest for 256-bit; the Strassen matrix multiplication is superior to the Ozaki scheme for 512-bit with $n \leq 1500$; and the Strassen matrix multiplication is the fastest for 768-bit computation.

From the above results, the 3M CGEMM method based on the Ozaki scheme with $d=14$ is the fastest for 256-bit calculations, whereas the Strassen matrix multiplication is the fastest for 768-bit calculations. Therefore, for matrices $A$ and $B$, the Ozaki scheme is the fastest in the precision range of 256 bits or less, and the Strassen matrix multiplication is the fastest for a precision of 768 bits or more.

%-----------------------------------%
% 
%-----------------------------------%
\subsection{Complex LU decomposition}

If the Ozaki scheme is faster than other matrix multiplication algorithms, it is expected to be effective for complex LU decomposition using matrix multiplication. Below we present the results of a benchmark test of complex LU decomposition using complex matrix multiplication based on the fast 3M method.

The corresponding $n$-dimensional linear system of equations is: 
\begin{equation}
	A\mathbf{x} = \mathbf{b}, \label{eqn:linear_eq}
\end{equation}
where $A\in\mathbb{C}^{n\times n}$, $\mathbf{x}\in\mathbb{C}^n$, and $\mathbf{b}\in\mathbb{C}^n$.

In this problem, we use $n=1024$. The elements of the complex coefficient matrix $A$ are given as random numbers with real and imaginary parts following a normal distribution in $[-1, 1]$. The exact solution of $\mathbf{x}$ is $(x)_k = k + k \mathrm{i}$ and $\mathbf{b} := A\mathbf{x}$ is a constant vector obtained with 2048-bit precision arithmetic.

\begin{figure}%[htb]
    \begin{center}
        \includegraphics[width=.35\textwidth]{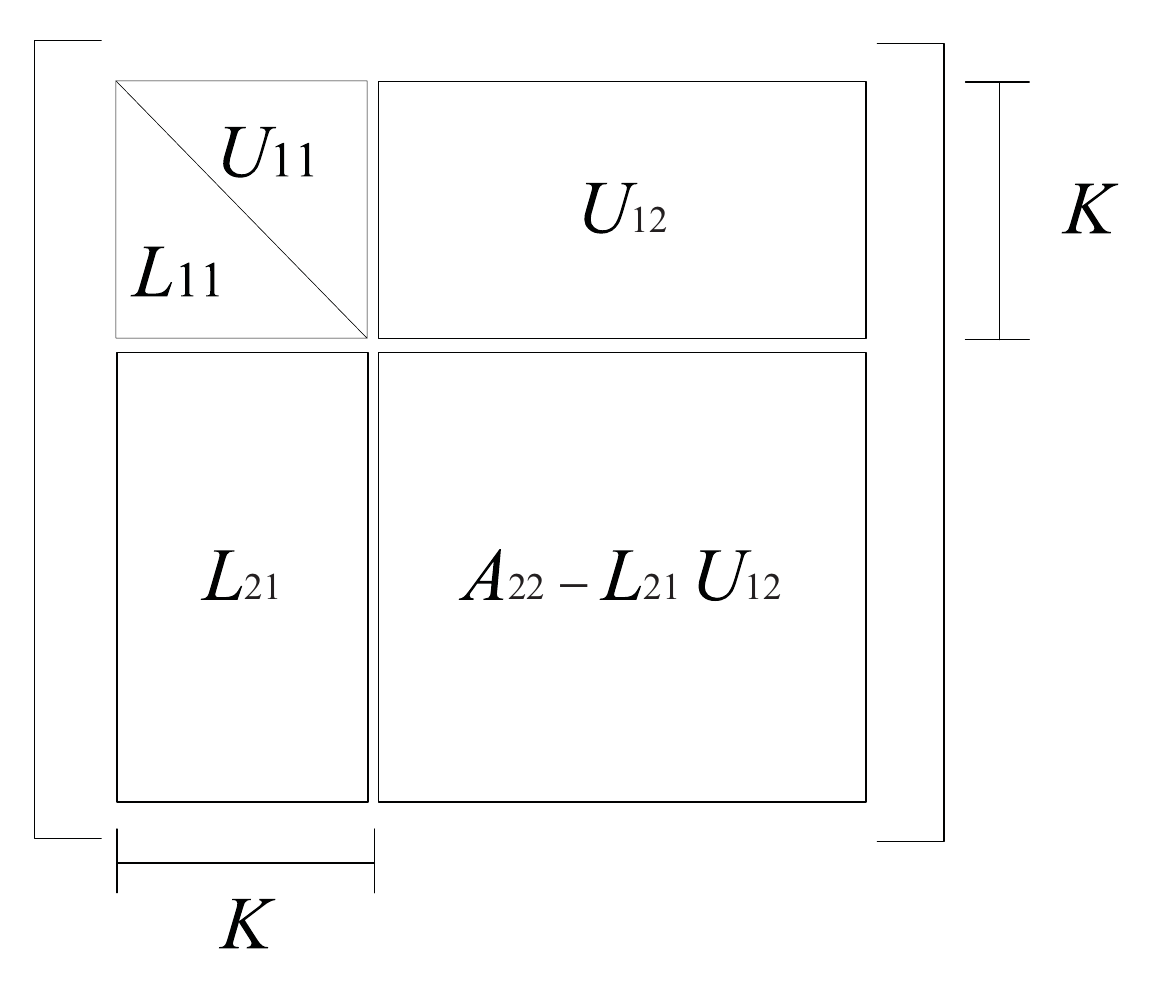}
        \caption{LU decomposition using matrix multiplication}
        \label{fig:lu}
    \end{center}
\end{figure}

Assuming that LU decomposition in the current LAPACK standard allows the use of fast matrix multiplication, a constant width $K$ is predefined as in \figurename\ref{fig:lu}, and the rectangular component $A - L_{21}U_{12}$ is updated for each $K$ column. Therefore, the complexity of matrix multiplication changes with respect to $K$; MPLAPACK's LU decomposition (Cgetrf) is also implemented using Cgemm in MPBLAS.

For fast LU decomposition, $L_{21}U_{12}$ must be fast in xGEMM.%; in DGEMM, the best $K$ must be specified for extremely small numbers of rows and columns, as performance deteriorates for these cases.

\begin{figure}%[htb]
    \begin{center}
        \includegraphics[width=.225\textwidth]{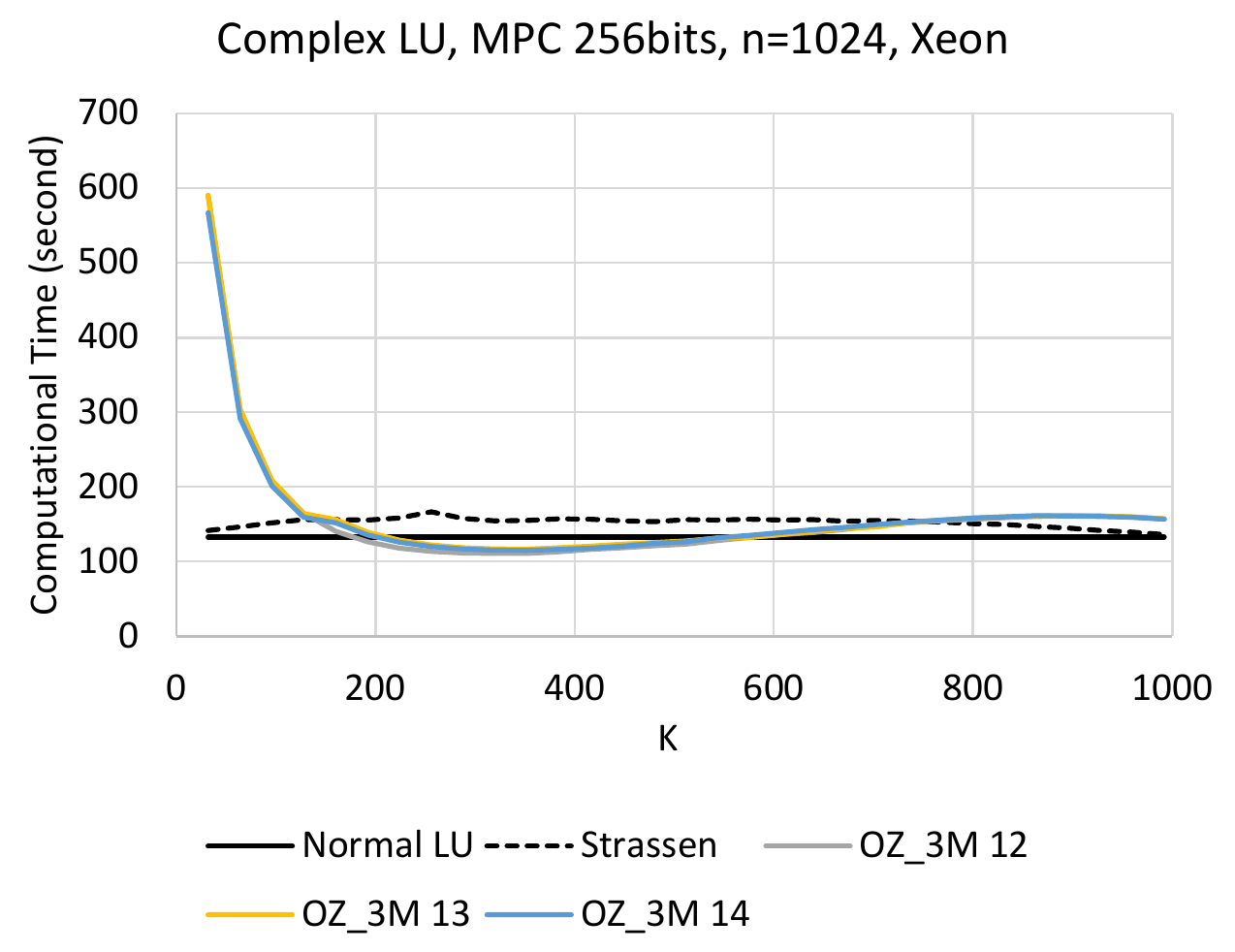}
        \includegraphics[width=.225\textwidth]{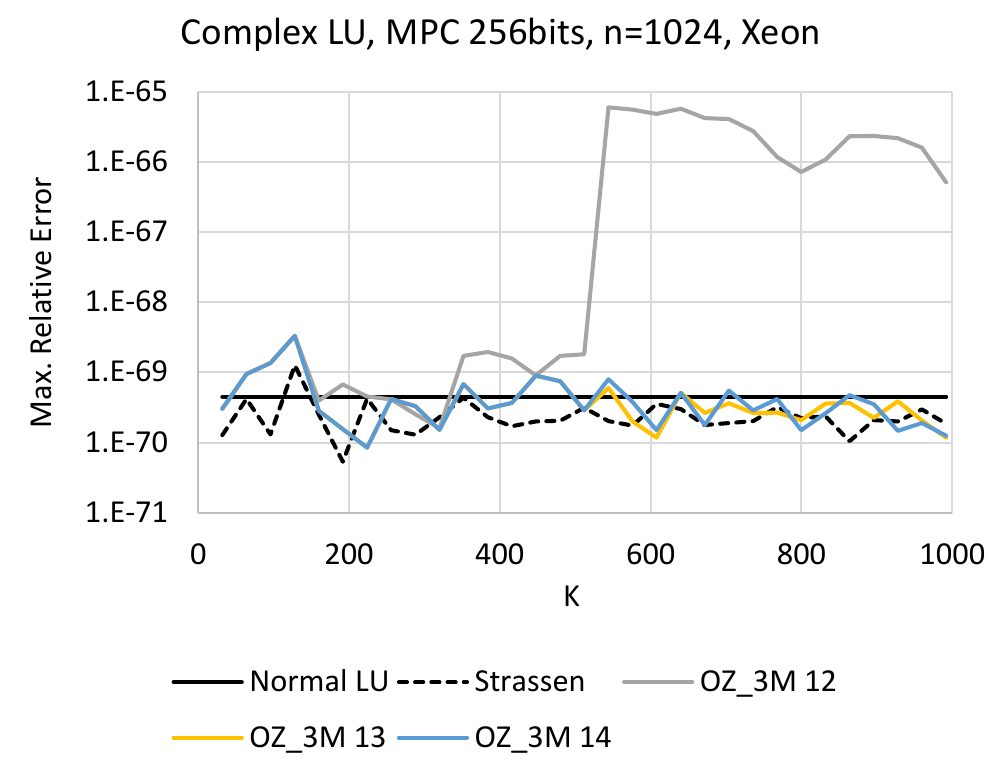}

		\includegraphics[width=.225\textwidth]{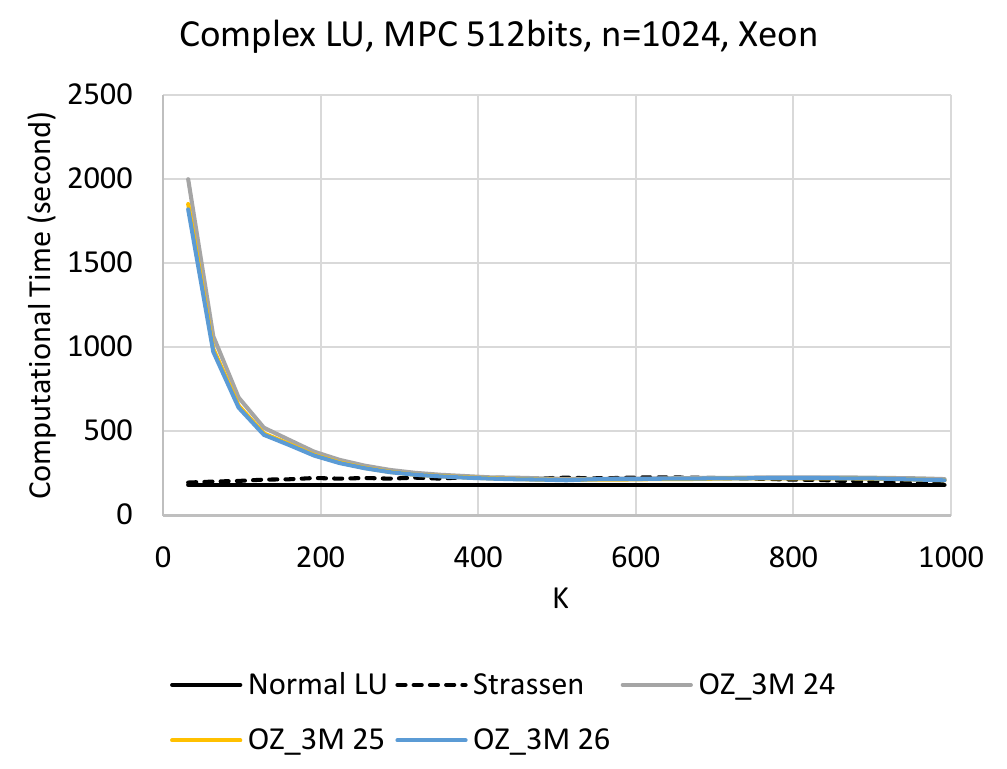}
        \includegraphics[width=.225\textwidth]{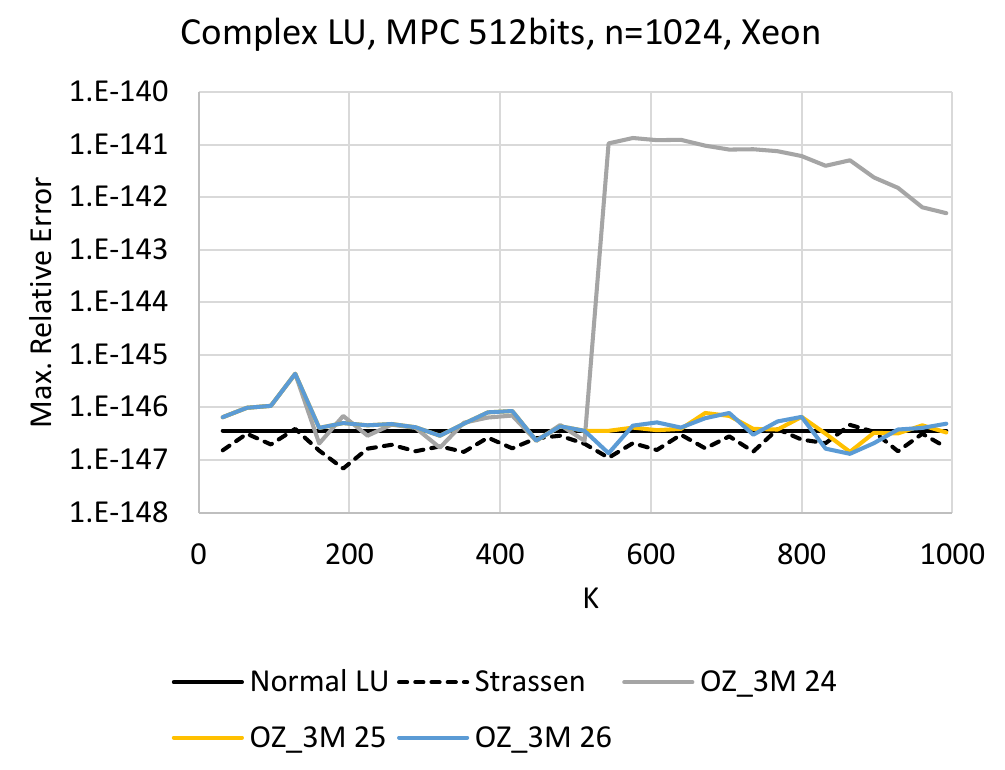}

		\includegraphics[width=.225\textwidth]{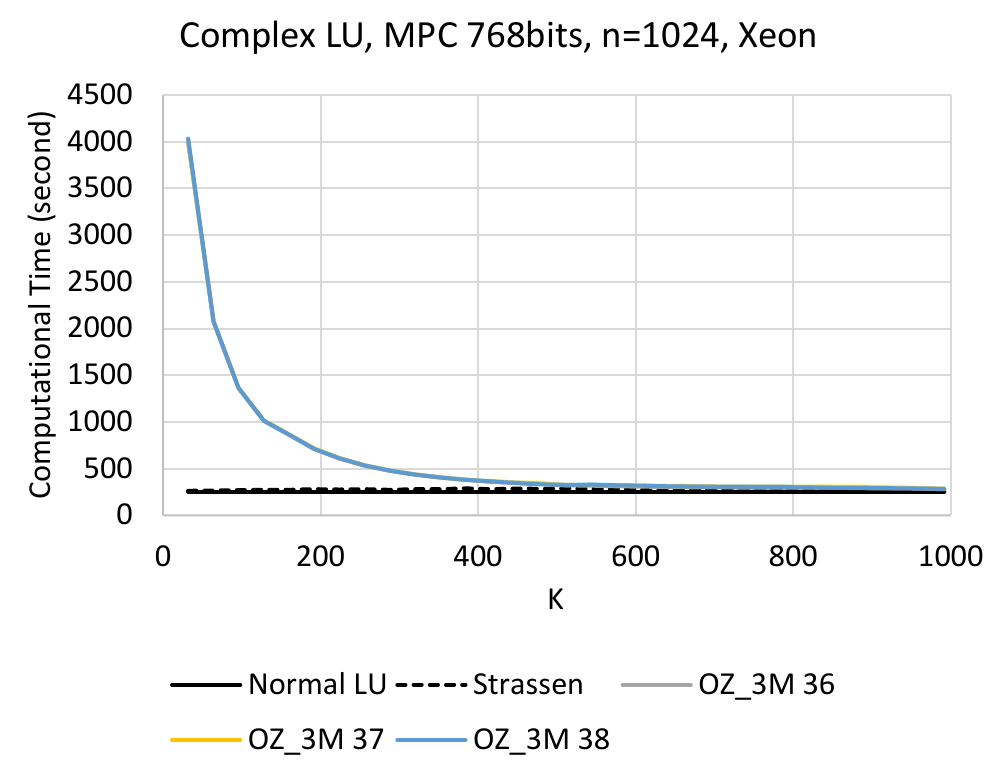}
        \includegraphics[width=.225\textwidth]{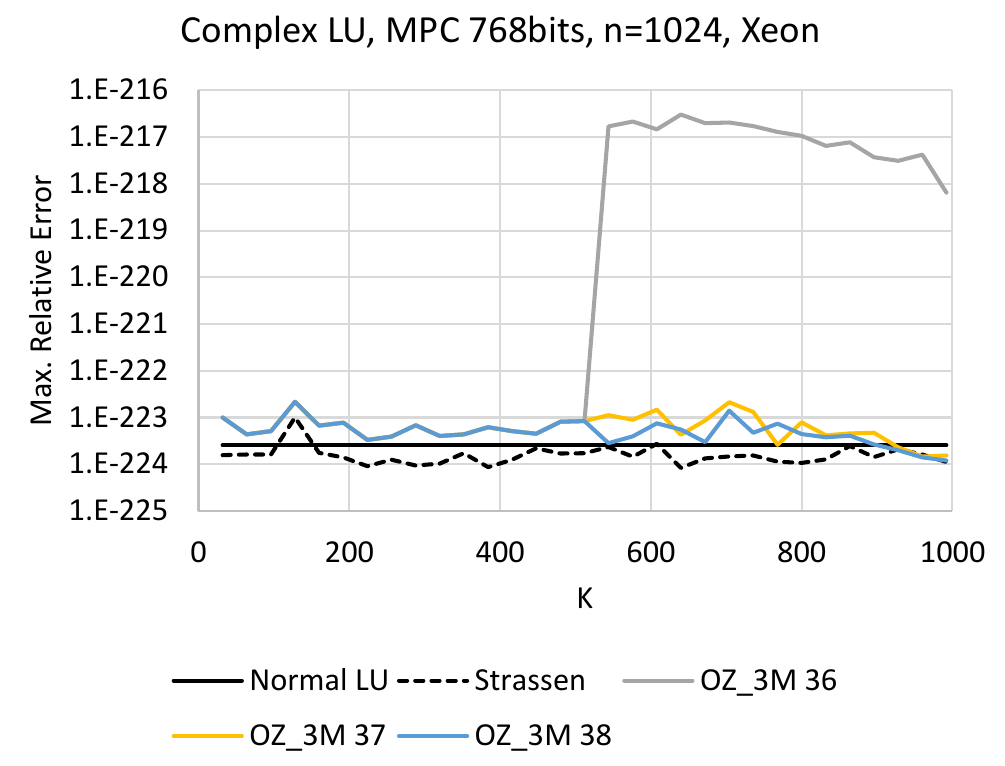}

		\caption{Computational time (left) and relative error (right) of complex LU decomposition on Xeon}
        \label{fig:clu_xeon}
    \end{center}
\end{figure}

\begin{figure}%[htb]
    \begin{center}
        \includegraphics[width=.225\textwidth]{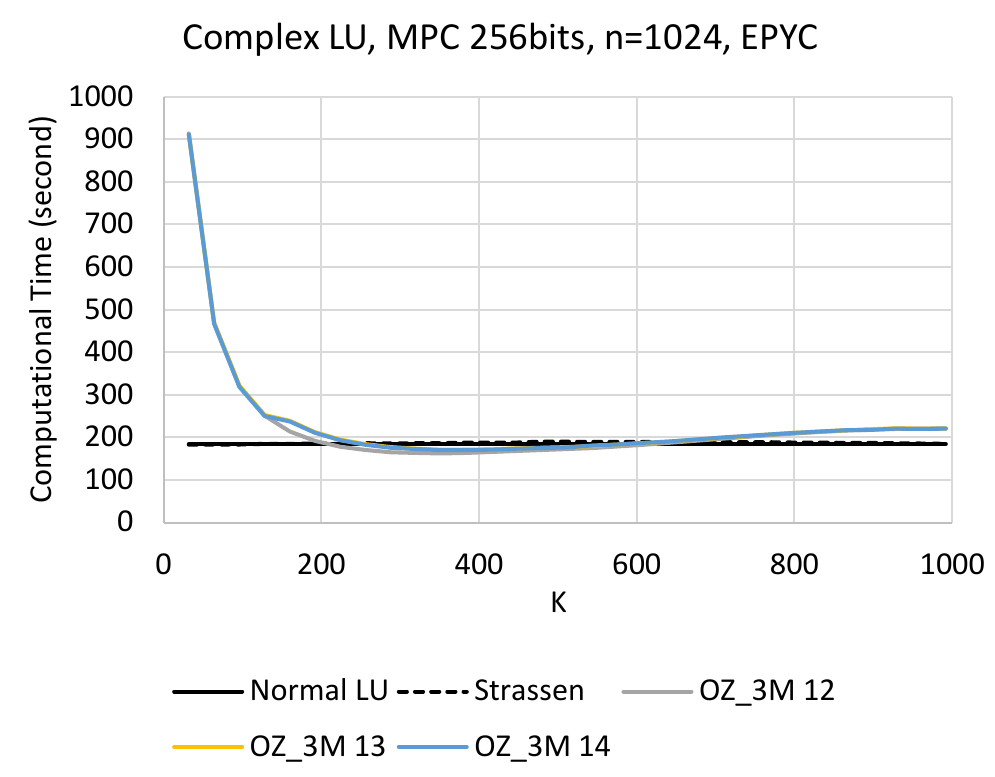}
        \includegraphics[width=.225\textwidth]{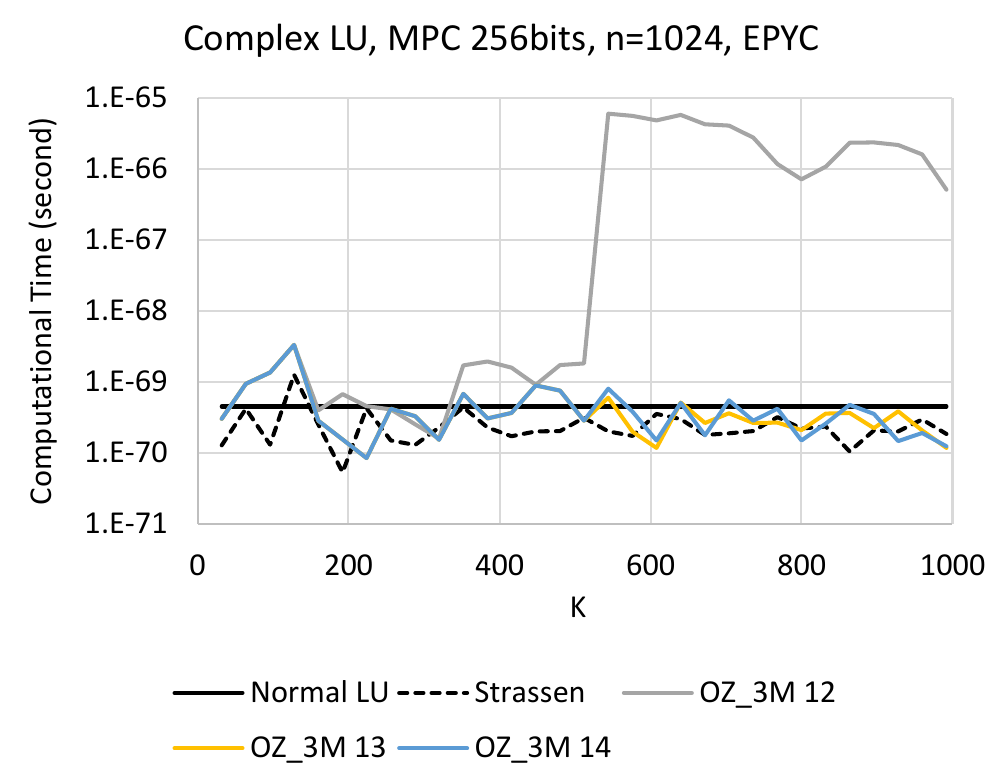}

		\includegraphics[width=.225\textwidth]{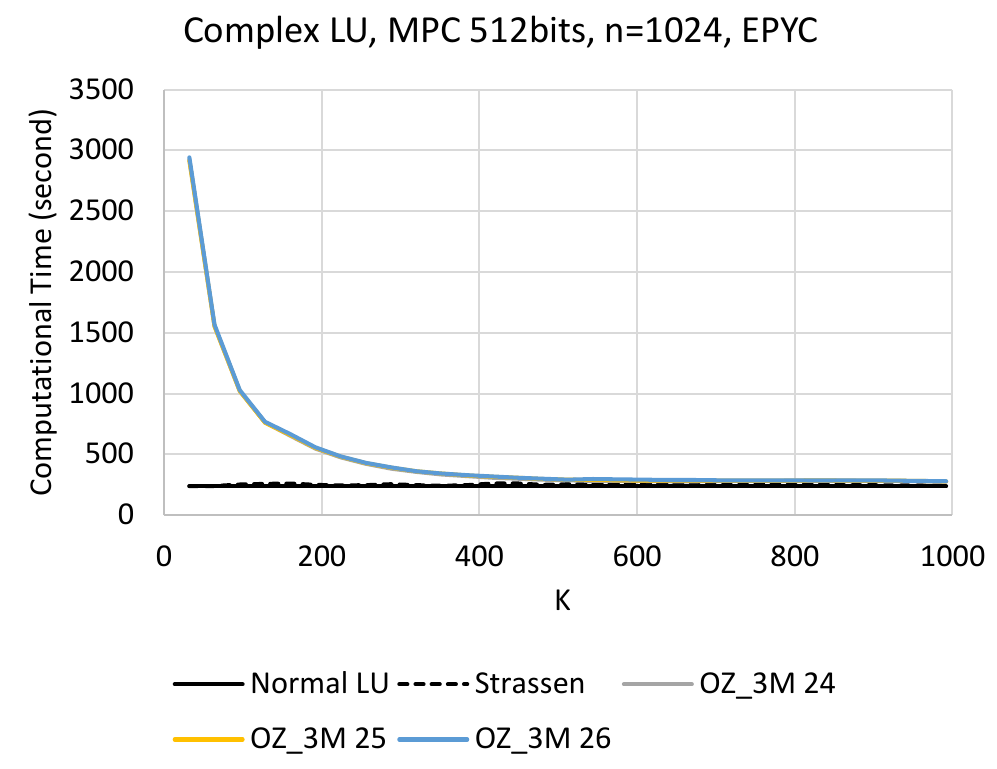}
        \includegraphics[width=.225\textwidth]{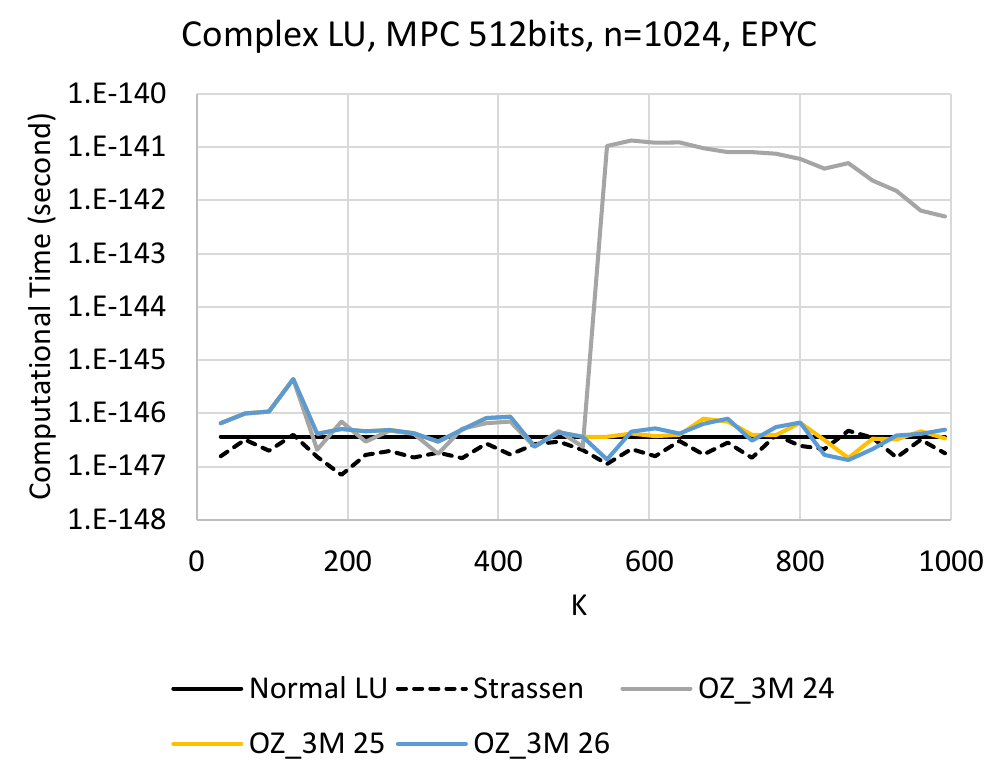}

		\includegraphics[width=.225\textwidth]{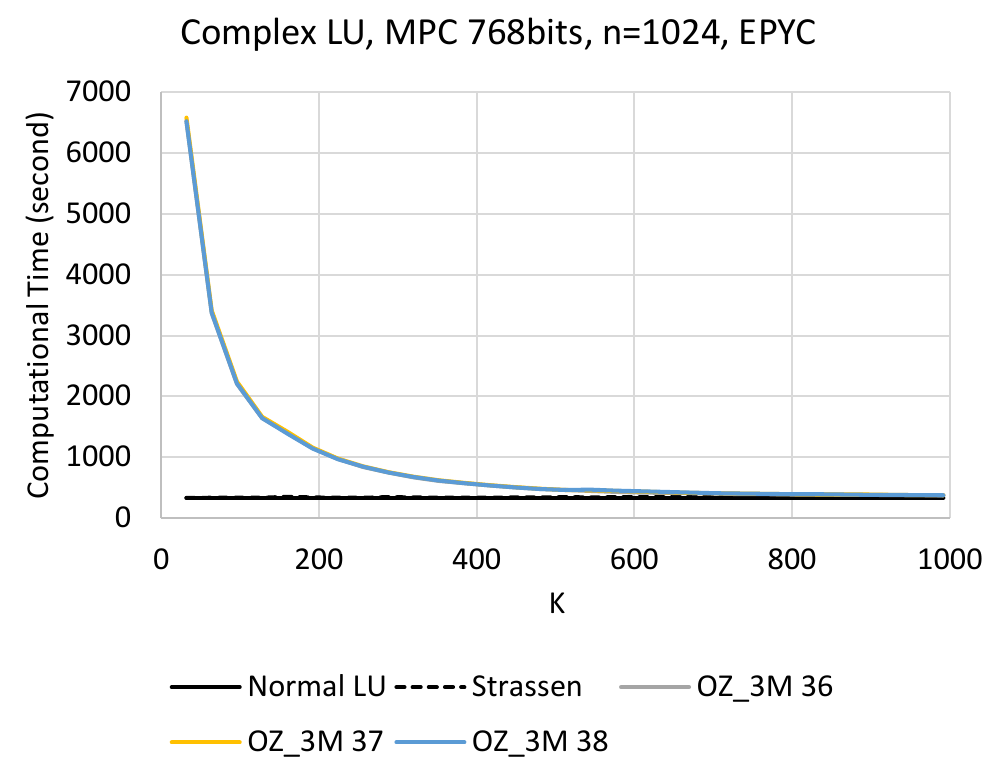}
        \includegraphics[width=.225\textwidth]{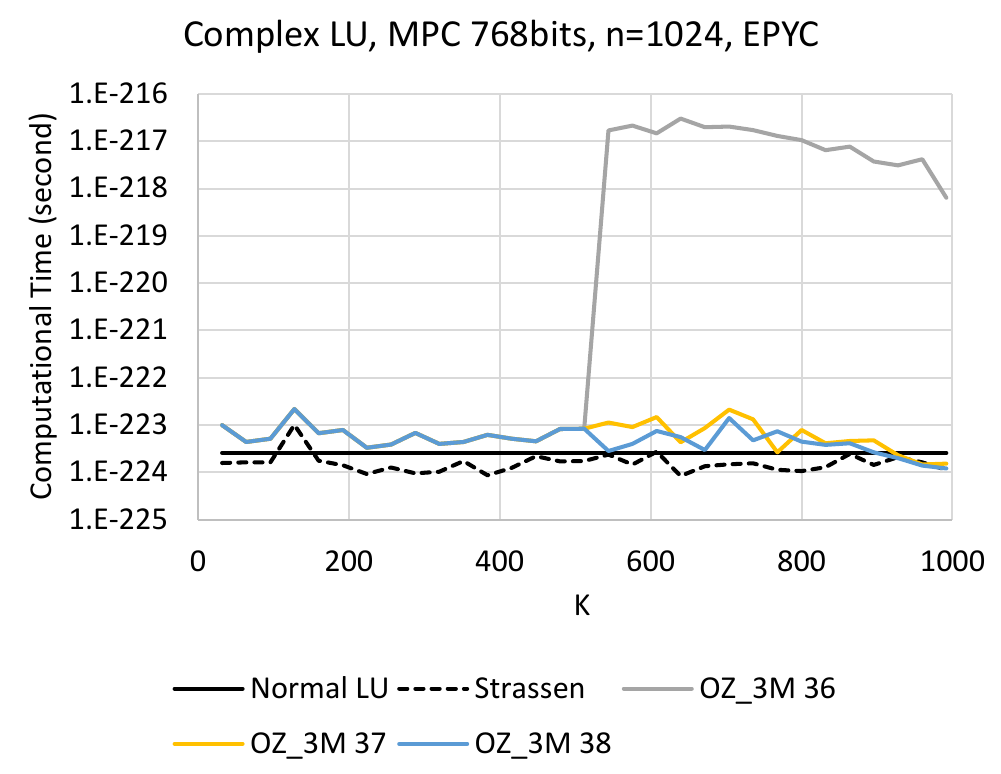}

		\caption{Computational time (left) and relative error (right) of complex LU decomposition on EPYC}
        \label{fig:clu_epyc}
    \end{center}
\end{figure}

In our benchmark test, the computational time and the maximum error of the numerical solution of the LU decomposition obtained by varying $K$ are measured using 256-, 512-, and 768-bit precision. The graphs of computational time (left) and relative error (right) of the numerical solution are displayed in \figurename\ref{fig:clu_xeon} and \figurename\ref{fig:clu_epyc}, respectively and are obtained by forward and backward substitutions. In addition, the results of normal column-wise LU decomposition (Normal LU, $K=1$) are also shown in these graphs.

There is no difference in the variation of the maximum relative error of the numerical solution when $K$ is changed in either the Xeon or EPYC environment. The trend in computational time is also approximately the same.

This leads us to conclude that,
\begin{itemize}
	\item Regardless of precision, the computational time is never less than that of Normal LU as long as the Strassen matrix multiplication is used. The only exception is 768-bit precision of LU decomposition on EPYC; however, the difference in computational time is extremely small.
	\item Using the Ozaki scheme, a 256-bit calculation can be performed in less computational time than for Normal LU with relatively small $K$. The number of divisions $d$ for which the relative error is minimal regardless of $K$ is $d=13$ for 256 bits, $d=25$ for 512 bits, and $d=37$ for 768 bits.
	\item Normal LU is the fastest for 512-bit and 768-bit calculations.
\end{itemize}

The $K$ and the maximum relative error associated with minimum computational time for each algorithm for each computing environment and precision are shown in \tablename \ \ref{table:min_comptime_mpfr_xeon} and \tablename\ \ref{table:min_comptime_mpfr_epyc}. The results of MPLAPACK's Cgetrf function are also shown for comparison.

\begin{table}
    \begin{center}
        \caption{Minimum computational time (s) of MPC LU decomposition and maximum relative errors of $\mathbf{x}$ on Xeon}\label{table:min_comptime_mpfr_xeon}
        \begin{tabular}{|c|c|c|c|c|}\hline
        Prec.    & Method        & $K$   & Second & Rel.Err. \\ \hline
        MPC      & Cgetrf        & N/A   & 188.1  & 2.9$\mbox{E}-71$\\
        256 bits & Normal LU     & 1     & 134.1  & 4.5$\mbox{E}-70$ \\
                 & Strassen      & 992   & 136.1  & 1.9$\mbox{E}-70$ \\
                 & OZ 13         & 352   & \underline{116.6}  & 6.9$\mbox{E}-70$ \\ \hline
        MPC      & Cgetrf        & N/A   & 236.7  & 6.4$\mbox{E}-149$ \\
        512 bits & Normal LU     & 1     & \underline{181.7}  & 3.7$\mbox{E}-147$ \\
                 & Strassen      & 992   & 187.6  & 8.0$\mbox{E}-147$ \\
                 & OZ 25         & 992   & 209.0  & 3.4$\mbox{E}-147$ \\ \hline
        MPC      & Cgetrf        & N/A   & 313.7  & 1.6$\mbox{E}-225$ \\
        768 bits & Normal LU     & 1     & \underline{253.0}  & 2.6$\mbox{E}-224$ \\
                 & Strassen      & 992   & 257.5  & 1.1$\mbox{E}-224$ \\
                 & OZ 37         & 992   & 287.6  & 1.5$\mbox{E}-224$ \\ \hline
        \end{tabular}
    \end{center}
\end{table}

\begin{table}
    \begin{center}
        \caption{Minimum computational time (s) of MPC LU decomposition and maximum relative errors of $\mathbf{x}$ on EPYC}\label{table:min_comptime_mpfr_epyc}
        \begin{tabular}{|c|c|c|c|c|}\hline
        Prec.    & Method        & $K$   & Second & Rel.Err. \\ \hline
        MPC      & Cgetrf        & N/A   & 275.9  & 2.9$\mbox{E}-71$\\
        256 bits & Normal LU     & 1     & 184.1  & 4.5$\mbox{E}-70$ \\
                 & Strassen      & 32    & 182.8  & 1.3$\mbox{E}-70$ \\
                 & OZ 13         & 384   & \underline{171.1}  & 3.1$\mbox{E}-70$ \\ \hline
        MPC      & Cgetrf        & N/A   & 327.7  & 6.4$\mbox{E}-149$ \\
        512 bits & Normal LU     & 1     & \underline{237.8}  & 3.7$\mbox{E}-147$ \\
                 & Strassen      & 64    & 240.7  & 3.2$\mbox{E}-147$ \\
                 & OZ 25         & 992   & 279.6  & 3.4$\mbox{E}-147$ \\ \hline
        MPC      & Cgetrf        & N/A   & 428.4  & 1.6$\mbox{E}-225$ \\
        768 bits & Normal LU     & 1     & 334.9  & 2.6$\mbox{E}-224$ \\
                 & Strassen      & 32    & \underline{333.5} & 1.6$\mbox{E}-224$ \\
                 & OZ 37         & 992   & 377.5  & $1.5\mbox{E}-224$ \\ \hline
        \end{tabular}
    \end{center}
\end{table}

The relative errors of Cgetrf, Normal LU, and Strassen, which use MPC only, are consistent across environments. Overall, the Cgetrf function is slow, but this may be due to the overhead of the C++ mpcomplex class and the performance of Cgemm in MPBLAS. These results indicate that computation speed can be further increased by using the Ozaki scheme and specifying the optimal $K$ and number of divisions $d$, for a precision of at least up to 256 bits.

%-----------------------------------%
% 
%-----------------------------------%
\section{Conclusion and future work}

It is demonstrated that optimization of complex matrix multiplication is possible on both Xeon and EPYC environments and that the Strassen matrix multiplication and implementation using the Ozaki scheme contribute to speed-up, especially compared to MPLAPACK. We also find that the 3M CGEMM method is fast and that the difference in accuracy is not noticeable for the problems examined in this study. Hence, there is no need to choose the 4M method, especially for multiple precision calculations. The use of fast complex matrix multiplication can also reduce the computational time of complex LU decomposition with the Ozaki scheme which is sufficiently faster than the Strassen matrix multiplication.

Future plans to expand this study include:
\begin{enumerate}
	\item Pursuing speed-up of arbitrary precision complex matrix multiplication. If the Strassen matrix multiplication is faster than the Ozaki scheme, a 3M CGEMM method using the Strassen algorithm is expected to be faster. Therefore, we plan to implement such a method and also parallelize using OpenMP. %By implementing these methods, further speed-up of complex LU decomposition can be achieved.
	\item Implementing and evaluating the performance of complex matrix multiplication in multi-component fixed-precision (double-double, triple-double or quadruple-double) arithmetic. In this case, the Ozaki scheme is expected to be advantageous in many cases with its relatively small computational accuracy.
	\item Investigating cases where the Ozaki scheme is useful for a wider range of problems. Complex matrix multiplication requires changing the number of divisions according to the difference in the order of absolute values of the real and imaginary parts of matrices. Therefore, we would like to investigate the precision and the number of divisions for real matrices by the Ozaki scheme, which are faster than the Strassen matrix multiplication. A method can then be established to set the optimal number of divisions for complex matrix multiplication based on these results.
\end{enumerate}

%---------------------------------------%
% Acknowledgement
%---------------------------------------%
\section*{Acknowledgment}
This work was supported by JSPS KAKENHI Grant Number 23K11127.

%---------------------------------------%
% Reference
%---------------------------------------%
%\bibliographystyle{IEEEtran}
%\bibliography{IEEEabrv,tkouya_utf8}

% Generated by IEEEtran.bst, version: 1.12 (2007/01/11)

\end{document}